\documentclass[a4paper,12pt]{article}

\usepackage{graphicx}
\usepackage{amssymb}

\def\NN{{\mathchoice {\mathrm{I \hspace{-0.2em} N}} 
    {\mathrm{I \hspace{-0.2em} N}} {\mathrm{I \hspace{-0.14em} N}} 
    {\mathrm{I \hspace{-0.14em} N}}}}
\def\ZZ{{\mathchoice {\mathsf{Z \hspace{-0.45em} Z}} {\mathsf{Z 
        \hspace{-0.45em} Z}} {\mathsf{Z \hspace{-0.32em} Z}} 
    {\mathsf{Z \hspace{-0.23em} Z}}}}
\def\QQ{{\mathchoice {\mathsf{I \hspace{-0.41em} Q}} {\mathsf{I
      \hspace{-0.47em} Q}} {\mathsf{I \hspace{-0.25em} Q}} 
      {\mathsf{I \hspace{-0.2em} Q}}}}

\newfam\gotfam
\font\twlgot=eufm10 at 12pt
\font\tengot=eufm10

\font\sevengot=eufm7
\textfont\gotfam=\twlgot
\scriptfont\gotfam=\tengot 
\scriptscriptfont\gotfam=\sevengot

\newtheorem{prop}  {Proposition}
\newtheorem{lemma}  {Lemma}
\newtheorem{theor}   {Theorem}
\newtheorem{truc} {Definition}

\newcommand{\be}  {\begin{equation}}
\newcommand{\ee}  {\end{equation}}
\newcommand{\bea} {\begin{eqnarray}}
\newcommand{\eea} {\end{eqnarray}}
\newcommand{\lp}  {\left(}
\newcommand{\rp}  {\right)}

\newcommand{\Br}  {\overline}

\newcommand{\cC}  {{\cal C}}

\newcommand{\cY}  {{\cal Y}}
\newcommand{\cO}  {{\cal O}}

\newcommand{\cI}  {{\cal I}}

\newcommand{\cF}  {{\cal F}}
\newcommand{\cR}  {{\cal R}}
\newcommand{\cA}  {{\cal A}}

\newcommand{\Om}  {\Omega}

\newcommand{\ta}  {\tau}
\newcommand{\ep}  {\epsilon}
\newcommand{\si}  {\sigma}

\newcommand{\Si}  {\Sigma}
\newcommand{\al}  {\alpha}

\newcommand{\La}  {\Lambda}

\newcommand{\ch}  {\chi}

\newcommand{\til}  {\tilde}

\def\Br{\overline}

\newcommand{\eqdef} {\stackrel{\rm def}{=}}

\def\endproof{\hfill\vrule height .6em width .6em depth
  0pt\goodbreak\vskip.25in}

\title{Grassmann-Berezin Calculus and Theorems of the 
Matrix-Tree Type }

\author{Abdelmalek Abdesselam \\
\\
{\small LAGA, Institut Galil\'ee, CNRS UMR 7539}\\
{\small  Universit{\'e} Paris XIII}\\
{\small Avenue J.B. Cl{\'e}ment, F93430 Villetaneuse, France}\\
{\small email: abdessel@math.univ-paris13.fr}
}

\begin{document}

\maketitle

{\abstract{
We prove two generalizations of the matrix-tree theorem.
The first one, a result essentially due to Moon for which we provide
a new proof, extends the ``all minors'' matrix-tree theorem
to the ``massive'' case where no condition on row or column sums  
is imposed. The second generalization,
which is new, extends the recently discovered
Pfaffian-tree theorem of Masbaum and Vaintrob into
a ``Hyperpfaffian-cactus'' theorem.
Our methods are noninductive, explicit and make critical use
of Grassmann-Berezin calculus that was developed for the needs of
modern theoretical physics.
}
}


\section{Introduction}
The matrix-tree theorem~\cite{Kirchhoff,Sylvester,Borchardt,Tutte}
is one of the most fundamental tools of
combinatorial theory.
Its applications are many, ranging from
electrical networks~\cite{Chen} to questions related to the partition function
of the Potts model in statistical mechanics~\cite{Sokal},
or to a recent conjecture
of Kontsevich regarding the number of points of varieties
defined by Kirchhoff spanning tree polynomials over finite
fields~\cite{Kontsevich,Stanley,Stembridge,Chung,Belkale}.
In its simplest instance, i.e. the {\em classical} matrix-tree theorem,
it says that the principal minors of a graph Laplacian
enumerate the spanning trees of the graph.
The matrix-tree theorem has many generalizations like
the ``all minors'' version~\cite{Chen,Chaiken,Moon}
and, more recently, the remarkable
Pfaffian-tree theorem of Masbaum and
Vaintrob~\cite{Masbaum1} whose motivation was the study of lowest
degree terms of Alexander-Conway polynomials of links and their
relation to Milnor invariants~\cite{Milnor,Masbaum2}.
We will prove both these generalizations of the matrix-tree theorem
using, in a critical and, we hope, illuminating manner,
what we call ``Grassmann-Berezin calculus'' in honor
of the main two inventors of the formalism.
This framework is also known as ``Fermionic integration''
or ``superanalysis''; see~\cite{Berezin,Deligne}
for mathematical precision,
or any modern textbook on quantum field theory for emphasis on computational
aspects.
We, by the way, would like to point out that the first example
of true Fermionic integration
(as opposed to mere determinant calculus) that we found in the
literature is the terrific letter~\cite{Clifford}
of Clifford to Sylvester,
where one can also find the ancestors of Feynman diagrams!
Grassmann-Berezin calculus is commonplace in modern
theoretical physics; it also strongly overlaps
with the more familiar exterior algebra.
We have nonetheless included, for the benefit of the reader,
a brief but self-contained review in Section 2,
where precise definitions are given and main properties are stated
without proof (see for instance ~\cite{Feldman}
or appendix B of~\cite{Salmhofer} for more detail).
In Section 3, we state and prove a generalization of the
all minors matrix-tree theorem for matrices that
are not necessarily symmetric with zero column sums, as is the case
for a graph Laplacian. Although not stated explicitly, this result is
essentially contained in~\cite{Moon} (see also~\cite{Chebotarev}).
Our proof is however a new one and serves as a
warm up session for Section 4,
where we provide a new generalization of the theorem of Masbaum
and Vaintrob, and express
a sum over spanning {\em cacti}, which is a hypergraph generalization
of the notion of tree
(our definition is different but related to the ones
in~\cite{Goulden,Elmaraki,Bona}), in terms of a Berezin integral
involving a collection of antisymmetric tensors which generalize
the ``matrix'' in ``matrix-tree''.
The mentioned Berezin integral, in a particular case
that includes the theorem of Masbaum and Vaintrob,
reduces to a Hyperpfaffian as considered, for instance,
in~\cite{Barvinok,Luque}.
The original proof~\cite{Masbaum1} of the Pfaffian-tree theorem used
an edge contraction induction.
Later, Hirschman and Reiner~\cite{Hirschman} found a noninductive
proof using a sign reversing involution (which, from the point of
view of combinatorial enumeration is more satisfactory).
Our proof, which is also noninductive and we hope even more
enlighting, builds on ideas by D. Brydges
related to the ``forest-root'' formula of~\cite{BrydgesI}.
The latter, is a promotion of an earlier formula of
Brydges and Wright~\cite{BrydgesW,Brydges}, 
which holds in a rather particular case,
into a much
more general ``fundamental theorem of calculus'', thereby illustrating
a  general principle
noticed in~\cite{Abdesselam} for similar identities.
We would like to add that the present paper is certainly not the last word
on possible generalizations of the matrix-tree theorem.
It seems, we dare say,
almost too easy to find more generalizations using the point of
view developed in this work, and we invite the reader to try
her/his own variation. A possible venue to explore is the generalization
of Theorem 2 below to cacti
that are not necessarily made of pieces with odd cardinality.
Another suggestion is to investigate what one could say for
tensors that are not
completely antisymmetric.
We believe that the best guide in trying to further extend
Theorem 2 is by having in mind a specific and relevant problem
from the theory of the symmetric group or
that of simplicial complexes.

\medskip
\noindent{\bf Acknowledgements :}
One of our motivations for the present work was to try to answer some
questions, related to the matrix-tree theorem and the $q\rightarrow 0$
limit of the Potts model,
raised by A. Sokal and generously submitted to our attention.
We thank D. Brydges for explaining to us the supersymmetric proof
of the ``forest-root'' formula of \cite{BrydgesI} and giving us a
good start by showing us how this proof translates
when applied to the case of the classical matrix-tree theorem.
We also thank G. Masbaum for his explanations as to the
knot-theoretical background of the Pfaffian-tree theorem.
Finally the support of the Centre National de la Recherche Scientifique
is most gratefully acknowledged.

\section{A review of Grassmann-Berezin calculus}
Let $R$ be a commutative ring with unit containing the 
field $\QQ$ of rational numbers.
Let $\ch_1,\ldots,\ch_n$ be a collection of letters.

\begin{truc}
The Grassmann algebra $R[\ch_1,\ldots,\ch_n]$, or simply $R[\ch]$,
is the quotient of the free noncommutative $R$-algebra
with generators $\ch_1,\ldots,\ch_n$, by the two-sided ideal
generated by the expressions
\be
\ch_i \ch_j+
\ch_j \ch_i
\ee
with $1\le i,j\le n$.
\end{truc}

In other words, the generators $\ch_i$ of $R[\ch]$
satisfy the {\em anticommutation} relations
\be
\ch_i \ch_j+
\ch_j \ch_i=0
\label{anti}
\ee
for all $i$ and $j$ in $[n]\eqdef \{1,\ldots,n\}$.
In particular, since $2$ is invertible,
one has $\ch_i^2=0$ for all $i\in [n]$.
The first important property of $R[\ch]$
is
\begin{prop}
$R[\ch]$ is a free $R$-module with basis given by the
$2^n$ monomials $\ch_{i_1}\ldots\ch_{i_p}$
with $0\le p\le n$, $1\le i_1<\cdots<i_p\le n$.
\end{prop}
A beautiful exercise we leave to the reader is to prove this
statement, which is the solution of a word problem, directly
from the definition, in a non inductive combinatorial way and
without using determinants, multilinear algebra or the universal
property that defines an exterior algebra.

As a result of the proposition any element $f\in R[\ch]$ can be
uniquely written as
\be
f=\sum_{p=0}^n
\sum_{1\le i_1<\cdots<i_p\le n}
f_{i_1\ldots i_p}
\ch_{i_1}\ldots\ch_{i_p}
\label{fexp}
\ee
with $f_{i_1\ldots i_p}\in R$.
One therefore has two natural gradings on the algebra
$ R[\ch]$ :
an $\NN$-grading by the number of factors $\ch$, i.e. the degree,
and a $\ZZ_2$-grading
$R[\ch]= R[\ch]_{\rm even} \oplus R[\ch]_{\rm odd}$
where $R[\ch]_{\rm even}$ (resp. $R[\ch]_{\rm odd}$) is generated, as
an $R$-module, by the monomials with an even (resp. odd) number of
factors.
A nonzero element $f$ which belongs to $R[\ch]_{\rm even}$
or $R[\ch]_{\rm odd}$ is said $\ZZ_2$-homogenous,
and its parity is $p(f)\eqdef 0$ in the first case
and $p(f)\eqdef 1$ in the second.
If $f$ , $g$ are  $\ZZ_2$-homogenous, one has
\be
f g=(-1)^{p(f)p(g)} g f
\ee
As a result one has the following most important
fact about Grassmann-Berezin calculus.

\begin{prop}
{\bf The Pauli exclusion principle :}
If $f$ is an odd element of $R[\ch]$, i.e. belongs to
$R[\ch]_{\rm odd}$, then
\be
f^2=0
\ee
\end{prop}
We will mostly use this property for $f$ homogenous of degree $1$,
where the physical terminology of `` Pauli exclusion principle''
most properly applies.
A consequence of the anticommutation relations (\ref{anti})
and the finiteness of the number $n$ of generators
is that every element of $R[\ch]_+$ (the set of elements with no
term in degree $0$) is nilpotent.
This allows, for instance, to define for any $f\in R[\ch]_+$
\be
exp(f)\eqdef
\sum_{p\ge 0}
\frac{1}{p!}f^p
\ee
since the series terminates after a finite number of terms.
We will however exclusively consider exponentials of {\em even}
elements, so that $e^{f+g}=e^f e^g$ holds.
For any $i$, $1\le i \le n$, we define the odd derivation
$\frac{\partial}{\partial\ch_i}$ acting to the right, as the
degree $-1$ $R$-linear map $R[\ch]\rightarrow R[\ch]$,
defined by the following action
on monomials $\ch_{i_1}\ldots\ch_{i_p}$
with $1\le i_1<\cdots<i_p\le n$.
We let
\be
\frac{\partial}{\partial\ch_i}
\ch_{i_1}\ldots\ch_{i_p}\eqdef 0
\ee
if $i\notin \{i_1,\ldots,i_p\}$ and
\be
\frac{\partial}{\partial\ch_i}
\ch_{i_1}\ldots\ch_{i_p}\eqdef 
(-1)^{\al-1}
\ch_{i_1}\ldots\ch_{i_{\al-1}}\ch_{i_{\al+1}}\ldots\ch_{i_p}
\ee
if there is an $\al$, $1\le\al \le p$ such that $i_\al=i$.

If $I=\{i_1,\ldots,i_p\}$, with $i_1<\cdots<i_p$
is  a subset of $[n]$, the Grassmann algebra
$R[\ch_I]\eqdef R[\ch_{i_1},\ldots,\ch_{i_p}]$ naturally embeds
into $R[\ch]=R[\ch_1,\ldots,\ch_n]$ and we will use
the corresponding identifications.
In particular, the degree zero part of $R[\ch]$
is identified with $R$.
As a result, for any injective map $\ta:[p]\rightarrow[n]$,
the $R$-linear composite map
$\frac{\partial}{\partial\ch_{\ta(1)}}
\circ\cdots\circ\frac{\partial}{\partial\ch_{\ta(p)}}$
can be viewed either as $R[\ch]\rightarrow R[\ch]$ or
$R[\ch]\rightarrow R[\ch_{I^{\rm c}}]$, where
$I^{\rm c}$ denotes the complement of $I\eqdef {\rm Im}\ \ta$
in $[n]$. Following F.A. Berezin, we use the {\em integral} notation
\be
\int {\rm d}\ch_{\ta(1)}\ldots {\rm d}\ch_{\ta(p)}\ f
\ee
for the image in $R[\ch_{I^{\rm c}}]$ of $f\in R[\ch]$ by the map
$\frac{\partial}{\partial\ch_{\ta(1)}}
\circ\cdots\circ\frac{\partial}{\partial\ch_{\ta(p)}}$.
Of particular importance is the case where $p=n$ and
$\ta(i)=n-i+1$ for any $i$, $1\le i\le p$.
If $f\in R[\ch]$ is written as in (\ref{fexp}) one then has
\be
\int {\rm d}\ch_{n}\ldots {\rm d}\ch_{1}\ f
=
f_{12\ldots n}
\ee
the ``top form'' coefficient of $f$.
Notice also that for any $f\in R[\ch]$ and
any permutation $\si$ of $[n]$,
\be
\int {\rm d}\ch_{\si(1)}\ldots {\rm d}\ch_{\si(n)}\ f
=
\ep(\si)
\int {\rm d}\ch_{1}\ldots {\rm d}\ch_{n}\ f
\ee
where $\ep(\si)$ denotes the signature of $\si$.
Now an easy consequence of the definitions is the following
\begin{prop}
If $n$ is an even integer and
$A$ is an $n\times n$ skew-symmetric matrix
with coefficients in $R$, and using the notation
$\ch A\ch\eqdef \sum_{i,j=1}^n \ch_i A_{ij}\ch_j$, one has
\be
\int {\rm d}\ch_{1}\ldots {\rm d}\ch_{n}\ e^{-\frac{1}{2}\ch A\ch}
= {\rm Pf}(A)
\ee
where ${\rm Pf}(A)$ denotes
the usual Pfaffian of $A$.
\end{prop}
We will also need
\begin{prop}
{\bf Fubini's theorem :}
Let $I=\{i_1,\ldots,i_p\}$ with $i_1<\cdots<i_p$
be a subset of $[n]$ and let $I^{\rm c}=\{j_1,\ldots,j_{n-p}\}$
with  $j_1<\cdots<j_{n-p}$,
then for any elements $f\in R[\ch_I]$ and
$g\in R[\ch_{I^{\rm c}}]$ we have, in the ring $R$,
the identity
\be
\int {\rm d}\ch_I
{\rm d}\ch_{I^{\rm c}}\ fg=
{(-1)}^{p(n-p)}
\lp
\int {\rm d}\ch_I\ f
\rp
\lp
\int {\rm d}\ch_{I^{\rm c}}\ g
\rp
\ee
where ${\rm d}\ch_I$ (resp. ${\rm d}\ch_{I^{\rm c}}$)
is shorthand for ${\rm d}\ch_{i_1}\ldots{\rm d}\ch_{i_p}$
(resp. ${\rm d}\ch_{j_1}\ldots{\rm d}\ch_{j_{n-p}}$ ).
\end{prop}

An important special case of the previous considerations is when
$n=2m$ is even and the variables $\ch_1,\ldots,\ch_n$ come in pairs
$\psi_i$, ${\Br \psi}_i$, $1\le i\le m$, i.e. when one works in the
Grassmann algebra $R[\psi,{\Br \psi}]\eqdef
R[\psi_1,\ldots,\psi_m,{\Br \psi}_1,\ldots,{\Br \psi}_m]$.
Although suggestive of complex conjugation, the bar is simply a
notation due to an extra combinatorial structure on the set $[n]$ that
labels the variables.
If $f\in R[\psi,{\Br \psi}]$, we introduce the notation
\be
\int ({\rm d}{\Br\psi} {\rm d}\psi)_{\rm ent}\ f\eqdef
\int
{\rm d}{\Br\psi}_1{\rm d}\psi_1
{\rm d}{\Br\psi}_2{\rm d}\psi_2
\ldots
{\rm d}{\Br\psi}_m{\rm d}\psi_m\ f
\ee
where ``ent'' is short for ``entangled form'' of the Berezin
integral of $f$.
The last result of Grassmann-Berezin calculus we need to recall
is the following.
\begin{prop}
If $A$ is any $m\times m$ matrix with coefficients
in $R$, and using the notation
${\Br\psi}A\psi\eqdef\sum_{i,j=1}^m {\Br\psi}_i A_{ij} \psi_j$, one
has
\be
\int ({\rm d}{\Br\psi} {\rm d}\psi)_{\rm ent}\ e^{-{\Br\psi}A\psi}=
{\rm det}(A)
\ee
More generally,
if $p$ is an integer $0\le p\le m$, and
$I=\{i_1,\ldots,i_p\}$,
$J=\{j_1,\ldots,j_p\}$ are two $p$-element subsets of $[m]$
where we made the choice of ordering $i_1<\cdots<i_p$ and
$j_1<\cdots<j_p$, if also $A_{I^{\rm c},J^{\rm c}}$
denotes the $(m-p)\times(m-p)$ matrix obtained
by erasing the rows of $A$ with index in $I$ and
the columns of $A$ with index in $J$, then
\be
\int ({\rm d}{\Br\psi} {\rm d}\psi)_{\rm ent}
\ (\psi_J {\Br \psi}_I)_{\rm ent}
e^{-{\Br\psi}A\psi}
=
(-1)^{\Si I+\Si J}{\rm det}(A_{I^{\rm c},J^{\rm c}})
\label{multicramer}
\ee
where $(\psi_J {\Br \psi}_I)_{\rm ent}\eqdef
\psi_{j_1}{\Br\psi}_{i_1}
\psi_{j_2}{\Br\psi}_{i_2}\ldots
\psi_{j_p}{\Br\psi}_{i_p}$,
$\Si I\eqdef i_1+\cdots+i_p$ and likewise for $\Si J$.
\end{prop}
Mind the inversion in the position of line and column variables.
Indeed, when $p=1$, $I=\{i\}$ and $J=\{j\}$, the quantity expressed
by either side of (\ref{multicramer}) is simply the matrix element
$({\rm com}\ A)_{ij}$ of the matrix of cofactors of $A$.
This allows, when $A$ is invertible, to elegantly rewrite
Cramer's rule as
\be
(A^{-1})_{ij}=
\frac{\int ({\rm d}{\Br\psi} {\rm d}\psi)_{\rm ent}\ \psi_i{\Br\psi}_j
e^{-{\Br\psi}A\psi}}
{\int ({\rm d}{\Br\psi} {\rm d}\psi)_{\rm ent}\ e^{-{\Br\psi}A\psi}}
\ee
in perfect analogy with the covariance of a complex Gaussian
probability measure.

\section{A generalization of the all minors matrix-tree theorem}
In this section we let $A=(A_{ij})_{1\le i,j\le n}$ be {\em any}
$n\times n$ matrix with entries in our ground ring $R$.
We will work in the Grassmann algebra
$R[\psi,{\Br \psi}]=
R[\psi_1,\ldots,\psi_n,{\Br \psi}_1,\ldots,{\Br \psi}_n]$.
Let $p$ be an integer, with $1\le p\le n$,
$I=\{i_1,\ldots,i_p\}$ and $J=\{j_1,\ldots,j_p\}$ be two $p$-element
subsets of $[n]$, fixed throughout this section, with
$i_1<\cdots<i_p$ and $j_1<\cdots<j_p$.
In the following a {\em forest} means a subset
of $K_n$ (the set of $2$-element subsets of $[n]$)
such that the associated graph, with vertex set $[n]$ and edge set
given by the forest itself, contains no cycle.
A {\em directed forest} $\cF$ is a set of pairs
$(u,v)\in[n]\times[n]$,
with $u\neq v$, such that if $(u,v)$ belongs to it, then $(v,u)$
does not, and such that the set $\{ \{u,v\}|(u,v)\in\cF \}$ is a
forest (undirected).
An edge $(u,v)$ in a directed forest $\cF$ is considered to be
oriented from $u$ to $v$. A directed forest $\cF$ (in fact its
associated undirected forest) naturally defines a partition $\Pi_\cF$
of $[n]$ into connected components.
$\cF$ restricts inside each block of $\Pi_\cF$ to a directed tree
that spans the block. With respect to the two sets $I$ and $J$,
a directed forest $\cF$ is called {\em admissible} if it satisfies the
following conditions :

- For any block $C\in\Pi_\cF$, either $C\cap (I\cup J)=\emptyset$
or both $C\cap I$ and $C\cap J$ are one-element sets.

- Inside any block $C\in\Pi_\cF$ that contains an element $i\in I$
and an element $j\in J$, all the edges of the corresponding directed
tree
are oriented {\em away} from $j$.

If $\cF$ is admissible, there is a unique permutation
$\si_\cF:[p]\rightarrow[p]$ such that for all $\al$,
$1\le \al\le p$, $j_\al$ and $i_{\si_\cF(\al)}$
are in the same component of $\Pi_\cF$.
The {\em signature} of $\cF$ is then defined as
$\ep(\cF)\eqdef \ep(\si_\cF)$.
Let $\cF$ be a subset of $[n]\times[n]$ and $\cR$ be a subset of
$[n]$. We say that the {\em pair} $(\cF,\cR)$ is {\em admissible}
if the following conditions are verified :

- $\cF$ is an admissible directed forest.

- Any $C\in\Pi_\cF$ which contains no element of
$I$ and $J$ has to contain a unique element of $\cR$.
Besides, $\cR$ has to be included in the union of such blocks
$C$.

- Inside any block $C$, like in the previous condition,
all the edges of the corresponding directed tree are oriented away
from the unique element of $C\cap \cR$ which plays the role
of a {\em root}.

Figure 1 shows an example of admissible pair
$(\cF,\cR)$. Here $n=16$, $I=\{3,7\}$, $J=\{2,8\}$,
$\cR=\{13,16\}$, and the directed forest is
\bea
\lefteqn{\cF=\{(2,4),(4,1),(4,7),(6,5),(6,3),(9,6),
} & & \nonumber\\
 & & (8,9),(9,10),(12,11),(13,12),(13,14),(13,15)\}
\eea
One also has
\bea
\lefteqn{\Pi_\cF=\left\{
\{1,2,4,7\},\{3,5,6,8,9,10\},
\right.
} & & \nonumber\\
 & & \left.
\{11,12,13,14,15\},\{16\}
\right\}
\eea
and $\ep(\cF)=-1$.

We can now state the following

\begin{figure}
\centering
\includegraphics[width=9.2cm,totalheight=6.5cm]{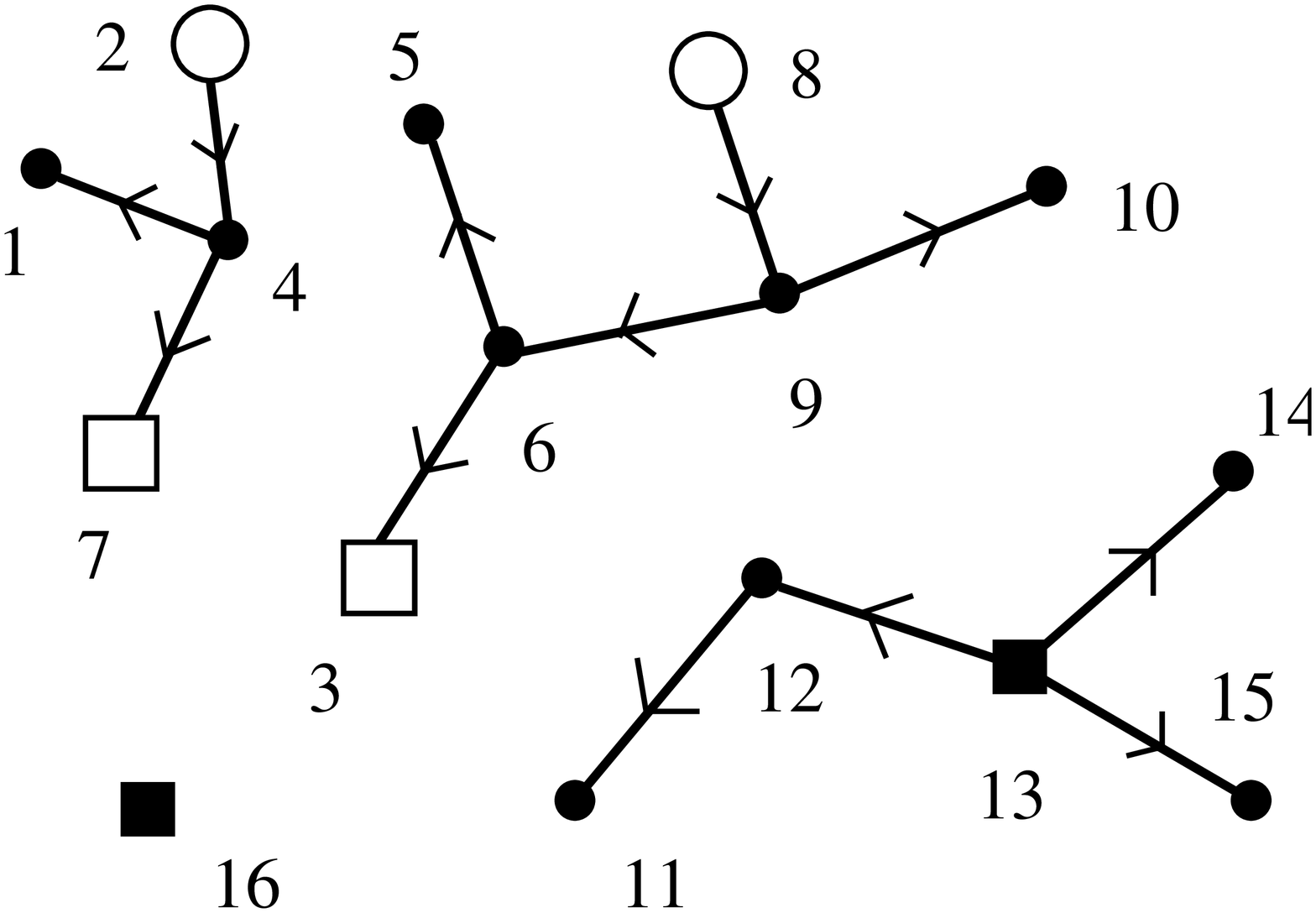}
\caption{An admissible pair $(\cF,\cR)$}
\end{figure}

\begin{theor}
\be
{\rm det}(A_{I^{\rm c},J^{\rm c}})
=
(-1)^{\Si I+\Si J}
\sum_{(\cF,\cR)\ {\rm admissible}}
\ep(\cF)
\prod_{j\in \cR}\lp
\sum_{i=1}^n
A_{ij}
\rp
\times
\prod_{(i,j)\in\cF}
\lp
-A_{ij}
\rp
\ee
\end{theor}

\noindent{\bf Proof :}
Let $\cI\eqdef
(-1)^{\Si I+\Si J}{\rm det}(A_{I^{\rm c},J^{\rm c}})$, which we
rewrite, thanks to Proposition 5, as
\be
\cI=
\int ({\rm d}{\Br\psi} {\rm d}\psi)_{\rm ent}
\ (\psi_J {\Br \psi}_I)_{\rm ent}
e^{-{\Br\psi}A\psi}
\ee
The trick, due to D. Brydges, that allows us to start is to write
\be
{\Br \psi}A\psi=
\sum_{j=1}^n
{\Br\psi}_j\lp
\sum_{i=1}^n
A_{ij}
\rp\psi_j
+
\sum_{i,j=1}^n
({\Br\psi}_i-{\Br\psi}_j)A_{ij}\psi_j
\ee
Let, for any $j$, $1\le j\le n$, $B_j\eqdef\sum_{i=1}^n A_{ij}$, one
then obtains
\be
\cI=
\int ({\rm d}{\Br\psi} {\rm d}\psi)_{\rm ent}
\ (\psi_J {\Br \psi}_I)_{\rm ent}
\exp\lp
-\sum_{j=1}^n B_j {\Br\psi}_j\psi_j
-\sum_{i,j=1}^n A_{ij}({\Br\psi}_i-{\Br\psi}_j)\psi_j
\rp
\ee
\be
=
\int ({\rm d}{\Br\psi} {\rm d}\psi)_{\rm ent}
\ (\psi_J {\Br \psi}_I)_{\rm ent}
\lp
\prod_{j=1}^n
e^{-B_j{\Br\psi}_j\psi_j}
\rp
\lp
\prod_{i,j=1}^n
e^{-A_{ij}({\Br\psi}_i-{\Br\psi}_j)\psi_j}
\rp
\ee
\be
=
\int ({\rm d}{\Br\psi} {\rm d}\psi)_{\rm ent}
\ (\psi_J {\Br \psi}_I)_{\rm ent}
\left[
\prod_{j=1}^n
\lp 1-B_j{\Br\psi}_j\psi_j \rp
\right]
\left[
\prod_{i,j=1}^n
\lp 1-A_{ij}({\Br\psi}_i-{\Br\psi}_j)\psi_j \rp
\right]
\ee
by the Pauli exclusion principle. We now expand to get
\be
\cI=
\sum_{(\cF,\cR)}
\lp
\prod_{j\in \cR} B_j
\rp
\lp
\prod_{(i,j)\in\cF}(-A_{ij})
\rp
\Om_{\cF,\cR}
\ee
where $\cF$ is {\em any} subset of $[n]\times[n]$, $R$ is
{\em any} subset of $[n]$ and we used the notation
\be
\Om_{\cF,\cR}\eqdef
\int ({\rm d}{\Br\psi} {\rm d}\psi)_{\rm ent}
\ (\psi_J {\Br \psi}_I)_{\rm ent}
\lp
\prod_{j\in\cR}
\left[ \psi_j {\Br\psi}_j \right]
\rp
\lp
\prod_{(i,j)\in\cF}
\left[ ({\Br\psi}_i-{\Br\psi}_j)\psi_j \right]
\rp
\ee
The theorem will now follow from the following
\begin{lemma}
$\Om_{\cF,\cR}=0$
unless the pair $(\cF,\cR)$
is admissible, in which case
$\Om_{\cF,\cR}=\ep(\cF)$.
\end{lemma}

\noindent{\bf Proof of the lemma :}
Trivially, if $(i,i)$ belongs to $\cF$, then the integrand of
$\Om_{\cF,\cR}$ contains a factor ${\Br\psi}_i-{\Br\psi}_i=0$ and
therefore $\Om_{\cF,\cR}$ vanishes.
Slightly less trivial is the fact that if both $(i,j)$ and $(j,i)$,
with $i\neq j$, belong to $\cF$ then again $\Om_{\cF,\cR}=0$.
Indeed, the integrand would then contain both the factors
$({\Br\psi}_i-{\Br\psi}_j)$ and $({\Br\psi}_j-{\Br\psi}_i)$
while $({\Br\psi}_i-{\Br\psi}_j)^2=0$
by the Pauli exclusion principle.
One more step down the ladder of triviality
takes us to the heart of the argument.
Suppose that the {\em undirected} graph associated to $\cF$ contains a
cycle, i.e. that for some $k\ge 3$ there is an injective map
$\ta:\ZZ/k\ZZ\rightarrow [n]$ such that for any $\al\in\ZZ/k\ZZ$,
$(\ta(\al),\ta(\al+1))$ or $(\ta(\al+1),\ta(\al))$ belongs to $\cF$.
Assume, for instance, that $(\ta(k),\ta(1))\in\cF$;
the alternate case can be treated in a similar vein.
Then, the integrand of $\Om_{\cF,\cR}$ contains the factor
\be
{\Br\psi}_{\ta(k)}-{\Br\psi}_{\ta(1)}=
({\Br\psi}_{\ta(k)}-{\Br\psi}_{\ta(k-1)})+\cdots+
({\Br\psi}_{\ta(2)}-{\Br\psi}_{\ta(1)})
\ee
Now, upon inserting this telescoping expansion of the factor
${\Br\psi}_{\ta(k)}-{\Br\psi}_{\ta(1)}$ into the integrand of 
$\Om_{\cF,\cR}$, the latter breaks into a sum of $(k-1)$ products.
For each of these products, there exists an $\al\in\ZZ/k\ZZ$
such that the factor $({\Br\psi}_{\ta(\al)}-{\Br\psi}_{\ta(\al-1)})$
appears {\em twice} : once with the $+$ sign from the telescopic
expansion of $({\Br\psi}_{\ta(k)}-{\Br\psi}_{\ta(1)})$, and once more
with a $+$ (resp. $-$) sign if $(\ta(\al),\ta(\al-1))$
(resp. $(\ta(\al-1),\ta(\al))$) belongs to $\cF$.
Again, the Pauli exclusion principle entails that $\Om_{\cF,\cR}=0$.

We now have reduced the discussion to the situation where $\cF$ is
a directed forest.
In this case, using Proposition 4, one can factor $\Om_{\cF,\cR}$
as $\Om_{\cF,\cR}=\ep\prod_{C\in\Pi_\cF} \Om_{\cF,\cR,C}$
where $\ep$ is a global sign we do not need to compute for the moment,
and for each $C\in\Pi_\cF$ of the form $C=\{c_1,\ldots,c_k\}$,
with $c_1<\cdots<c_k$,
\bea
\Om_{\cF,\cR,C} & \eqdef &
\int ({\rm d}{\Br\psi}_C {\rm d}\psi_C)_{\rm ent}
\lp
\prod_{j\in J\cap C} \psi_j
\rp
\lp
\prod_{i\in I\cap C} {\Br\psi}_i
\rp
\nonumber\\
 & & \lp
\prod_{j\in\cR\cap C}\lp \psi_j{\Br\psi}_j\rp
\rp
\lp
\prod_{(i,j)\in\cF_C}
\lp {\Br\psi}_i-{\Br\psi}_j\rp\psi_j
\rp
\label{component}
\eea
where any ordering of the factors in $\prod_{j\in J\cap C} \psi_j$
and $\prod_{i\in I\cap C} {\Br\psi}_i$
will do (eventual signs being absorbed in $\ep$),
$({\rm d}{\Br\psi}_C {\rm d}\psi_C)_{\rm ent}$ is shorthand
for
\[
{\rm d}{\Br\psi}_{c_1} {\rm d}\psi_{c_1}
{\rm d}{\Br\psi}_{c_2} {\rm d}\psi_{c_2}\ldots
{\rm d}{\Br\psi}_{c_k} {\rm d}\psi_{c_k}
\]
and
$\cF_C\eqdef \cF\cap(C\times C)$ is a spanning directed tree on the
vertex set $C$.
Note that, in order to have $\Om_{\cF,\cR,C}\neq 0$, there needs
to be exactly $k$ factors $\psi$ and as many factors ${\Br\psi}$
in the integrand. Since $\cF$ necessarily has $k-1$ edges,
the last product in (\ref{component}) already contributes
$k-1$ factors $\psi$ and $k-1$ factors ${\Br\psi}$.
This places severe restrictions on the sets $J\cap C$, $I\cap C$ and
$\cR\cap C$.
Either $J\cap C$ and $I\cap C$ are singletons and $\cR\cap
C=\emptyset$
in which case we say that $C$ is of {\em type I},
or $J\cap C=I\cap C=\emptyset$ and $\cR\cap C$ is a singleton
in which case we say that $C$ is of {\em type II}.
Note that the definition of $\Om_{\cF,\cR,C}$ is now unambiguous
since there is no problem of ordering the factors in
$\prod_{j\in J\cap C} \psi_j$
and $\prod_{i\in I\cap C} {\Br\psi}_i$ anymore.
One can readily check that the global sign $\ep$ is equal to the
signature $\ep(\cF)$ of $\cF$.
Finally we need to evaluate the expressions $\Om_{\cF,\cR,C}$ in the
two following cases.

\noindent{\bf 1st case : C of type I}

If $C\cap I=\{i\}$ and $C\cap J=\{j\}$ then
\be
\Om_{\cF,\cR,C}
=
\int ({\rm d}{\Br\psi}_C {\rm d}\psi_C)_{\rm ent}
\psi_j
{\Br\psi}_i
\lp
\prod_{(\al,\beta)\in\cF_C}\lp
{\Br\psi}_\al-{\Br\psi}_\beta
\rp \psi_\beta
\rp
\label{typeI}
\ee
First, note that there is a unique shortest path,
we call the {\em backbone},
joining $i$ and $j$ in the undirected tree associated to $\cF_C$.
Second, we need to inductively expand the product in (\ref{typeI})
starting from the leaves of the branches that hang from the backbone.
Let $\al\in C$ be such a leaf.
Then either $(\al,\beta)\in\cF_C$ or
$(\beta,\al)\in\cF_C$ for some $\beta\in C$.
In the first case, we write the corresponding factor
as $-\psi_\beta{\Br\psi}_\al+\psi_\beta{\Br\psi}_\beta$ and notice
that one cannot obtain the variable $\psi_\al$ in the integrand
and therefore $\Om_{\cF,\cR,C}=0$.
In the second case we get a factor
$-\psi_\al{\Br\psi}_\beta+\psi_\al{\Br\psi}_\al$.
If we keep the term $-\psi_\al{\Br\psi}_\beta$ in the expansion
then again there is no way of obtaining the factor ${\Br\psi}_\al$.
Therefore, to get a nonzero contribution, the edge containing
the leaf $\al$ has to be oriented {\em towards} $\al$ and we have
no choice but to select the term $\psi_\al{\Br\psi}_\al$
in the expansion.
Similarly to the Pr\"ufer coding of Cayley trees, we continue
this rewriting of $\Om_{\cF,\cR,C}$ by treating the
$({\Br\psi}_\al-{\Br\psi}_\beta)\psi_\beta$ factors corresponding to
the leaves, then to the vertices that become leaves after the first
generation leaves have been plucked out etc. until we arrive at the
backbone which plays the role of a root.
We then get $\Om_{\cF,\cR,C}=0$
unless all the edges, that are not on the backbone, are oriented away
from it, in which case
\be
\Om_{\cF,\cR,C}=
\int ({\rm d}{\Br\psi}_C {\rm d}\psi_C)_{\rm ent}
\lp
\prod_{\al\notin B} \psi_\al{\Br\psi}_\al
\rp
\La_B
\ee
where $B$ is the set of vertices on the backbone and $\La_B$ is an
expression to be defined as follows.
Let $k$ be an integer $k\ge 1$ and $\ta:[k]\rightarrow B$
be a bijective map such that $\ta(1)=j$ and
$\ta(k)=i$, and for any $l$, $1\le l\le k-1$,
$(\ta(l),\ta(l+1))$ or $(\ta(l+1),\ta(l))$
belongs to $\cF_C$.
If $(\ta(l),\ta(l+1)\in\cF_C$ we say that $l$ is {\em good}, and if
$(\ta(l+1),\ta(l)\in\cF_C$ we say that $l$ is {\em bad}.
Now
\bea
\lefteqn{
\La_B=\psi_j {\Br\psi}_i\times
} & & \nonumber\\
 & & \lp
\prod_{{1\le l\le k-1}\atop{l\ {\rm good}}}
\lp {\Br\psi}_{\ta(l)}-{\Br\psi}_{\ta(l+1)} \rp\psi_{\ta(l+1)}
\rp
\lp
\prod_{{1\le l\le k-1}\atop{l\ {\rm bad}}}
\lp {\Br\psi}_{\ta(l+1)}-{\Br\psi}_{\ta(l)} \rp\psi_{\ta(l)}
\rp
\eea
Let $l=1$, if $l$ is bad, then the corresponding factor is
$({\Br\psi}_{\ta(2)}-{\Br\psi}_j)\psi_j$.
Since $\psi_j^2=0$, we would then have $\Om_{\cF,\cR,C}=0$.
So $l$ has to be good and when we expand the corresponding
factor $({\Br\psi}_j-{\Br\psi}_{\ta(2)})\psi_{\ta(2)}=
{\Br\psi}_j\psi_{\ta(2)}-{\Br\psi}_{\ta(2)}\psi_{\ta(2)}$ we need
to keep the first term ${\Br\psi}_j\psi_{\ta(2)}$ otherwise
${\Br\psi}_j$ would not appear in the integrand and
$\Om_{\cF,\cR,C}$ would vanish. We then treat similarly
$l=2,3,\ldots,k-1$ to obtain that $\Om_{\cF,\cR,C}=0$ unless
also all the edges of the backbone are directed away from $j$,
in which case we are left with
\be
\Om_{\cF,\cR,C}=
\int ({\rm d}{\Br\psi}_C {\rm d}\psi_C)_{\rm ent}
\prod_{\al\in C}(\psi_\al{\Br\psi}_\al)
=1
\ee

\noindent{\bf 2nd case : C of type II}

It is exactly the same argument as in the previous case
in the degenerate situation where the backbone is reduced
to a single vertex $u$, with $u$ being the unique element of
$\cR\cap C$.

It is now simply a matter of checking our previous definitions of
admissibility to conclude the proof of the lemma.
\endproof

Now Theorem 1 follows immediately.
\endproof

\noindent{\bf Remark : }
The more familiar all minors matrix-tree theorem, as
one can find in~\cite{Chaiken}, is the ``massless'' particular case
of theorem 1 where the column sums of the matrix $A$
are zero, and where the only set of roots $\cR$ that gives
a nonzero contribution is $\cR=\emptyset$.
Theorem 1 was not explicitly stated in~\cite{Moon}; it however follows from
the general determinant expansion therein. Other related results
are reviewed in~\cite{Chebotarev}.

\section{A hyperpfaffian-cactus theorem :}
In this section we suppose $n$ is an odd positive integer,
and we work in the Grassmann algebra $R[\ch]=R[\ch_1,\ldots,\ch_n]$.
Suppose we are given for any {\em odd} integer $k$, $3\le k\le n$,
a completely antisymmetric tensor
$(y_{\al_1\ldots\al_k}^{[k]})_{(\al_1,\ldots,\al_k)\in[n]^k}$
with entries in the ground ring $R$.
It is simply a multidimensional analog of a matrix , and
complete antisymmetry means that for any
$(\al_1,\ldots,\al_k)\in[n]^k$ and any permutation $\si$ of the
set $[k]$
\be
y_{\al_{\si(1)}\ldots\al_{\si(k)}}^{[k]}=\ep(\si)
y_{\al_1\ldots\al_k}^{[k]}
\ee
where $\ep(\si)$ denotes the signature of $\si$.
Let ${\til \cO}_n$ denote the set of all subsets of $[n]$ which
have odd cardinality greater than or equal to $3$.
To any subset $\cA$ of ${\til \cO}_n$, we can associate an ordinary bipartite
graph $G(\cA)$ with vertex set partitioned into the disjoint union
of $\cA$ and $[n]$, and
edge set equal to the set of all pairs $(A,j)\in\cA\times [n]$
such that $j\in A$.
We say that $\cA$ is an {\em odd cactus} or simply a cactus,
if and only if $G(\cA)$ is a tree that connects all the
vertices of $[n]$.
Let ${\cO}_{n,k}$ denote the set of all sequences
$\al=(\al_1,\ldots,\al_k)$ of odd length $k$ made of {\em distinct}
elements of $[n]$, and let
\be
{\cO}_n\eqdef\cup_{{3\le k\le n}\atop{k\ {\rm odd}}}
{\cO}_{n,k}
\ee
To each sequence $\al=(\al_1,\ldots,\al_k)\in {\cO}_n$ we associate the
unordered set ${\til\al}\eqdef \{\al_1,\ldots,\al_k\}$ in ${\til \cO}_n$.
Let $\cC$ be a subset of ${\cO}_n$, we say that $\cC$ is a
{\em refined cactus} if and only if the following two conditions are
satisfied :

- For any distinct elements $\al$ and $\beta$ of $\cC$,
the sets ${\til\al}$ and ${\til\beta}$ are also distinct.

- $\cA(\cC)\eqdef \{ {\til\al}\in{\til \cO}_n | \al\in\cC \}$ is a cactus.

\begin{figure}
\centering
\includegraphics[width=9.2cm,totalheight=6.5cm]{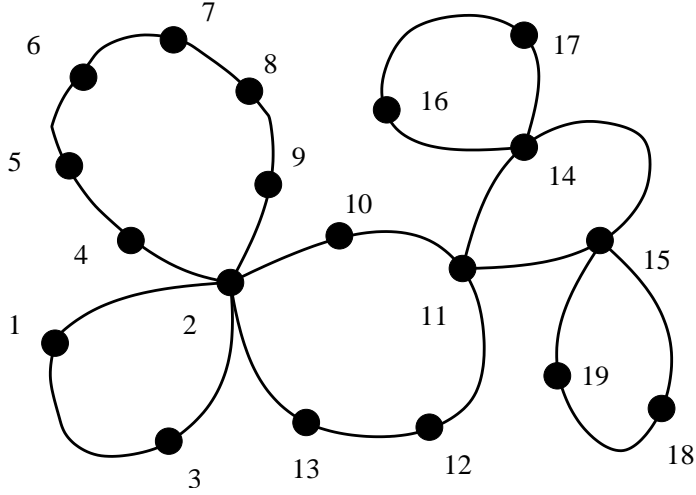}
\caption{A cactus}
\end{figure}

Figure 2 shows a possible representation of what we called an odd
cactus.
Here $n=19$, and the odd cactus is
\bea
\lefteqn{\cA=\left\{
\{1,2,3\},\{2,4,5,6,7,8,9\},\{2,10,11,12,13\},
\right.
} & & \nonumber\\
 & & \left.
\{11,14,15\},\{14,16,17\},\{15,18,19\}
\right\}
\eea
Note that $\cA$ is a set of unordered subsets of $[n]$.
For instance, if one would arbitrarily permute the
labels $4,5,6,7,8,9$ on the picture,
the odd cactus would still be the same.
In fact, the cyclic structure of the ``lobes'' of the cactus,
in the planar representation of Figure 2, is more relevant for what we
called a refined cactus.
For instance, a refined cactus $\cC$, corresponding to the
previous odd cactus $\cA$, and for which Figure 2 is a more
faithful representation is
\bea
\lefteqn{\cC=\left\{
(2,3,1),(6,7,8,9,2,4,5),(11,12,13,2,10),
\right.
} & & \nonumber\\
 & & \left.
(11,14,15),(17,14,16),(15,18,19)
\right\}
\eea
where the ordering of any sequence $\al\in\cC$
agrees with {\em clockwise} rotation on the corresponding
lobe of the cactus. Note that, even with this rule,
Figure 2 is still ambiguous in specifying
a refined cactus since one still has to chose the starting point
of every sequence $\al$.
Indeed, there is $3^4\times 5\times 7=2835$
possible refined cacti corresponding to Figure 2.

Let $i$ be a fixed vertex of $[n]$ which will play the role of a root
and let $\cC$ be a refined cactus.
For any $\al=(\al_1,\ldots,\al_k)$ in $\cC$,
there is a unique shortest path in the bipartite tree graph
$G(\cA(\cC))$ going from ${\til\al}$ to $i$.
The first vertex of $[n]$ one meets along this path starting from
${\til\al}$ is called the {\em local root} of $\al$ and is of the form
$\al_s$ for a unique index $s$, $1\le s\le k$.
We then define the {\em circulation } of $\al$ as the sequence
\be
{\hat\al}\eqdef (\al_{s+1},\al_{s+2},\ldots,\al_k,\al_1,\al_2,\ldots
\al_{s-1})
\ee
which has an even lenght $k-1$.
Now choose an ordering of $\cC$, and define by {\em concatenation}
a sequence $\pi$ by putting the $i$ first and then successively
all the sequences ${\hat\al}$, for $\al\in\cC$ according to the
chosen ordering of $\cC$.
Note that $\pi=(\pi_1,\ldots,\pi_n)$ is a permutation of the sequence
$(1,2,\ldots,n)$.
For example, if $\cC$ is a refined cactus represented by Figure 2
and if one chooses $i=10$ as a root,
then there is $6!=720$ possible sequences $\pi$ to choose from, one of which
is, for instance
\bea
\lefteqn{\pi=(10,11,12,13,2,3,1,4,5,6,
} & & \nonumber\\
 & & 7,8,9,18,19,14,15,16,17)
\eea

We let $\ep_{i,\cC}$ denote the signature of $\pi$.
This is well defined, since changing the ordering of $\cC$ amounts to
rigidly moving around the ${\hat\al}$'s which are all of {\em even}
length.
We can now define the {\em amplitude} of a refined cactus $\cC$,
with respect to the choice of root $i$ as
\be
\cY_{i,\cC}\eqdef \ep_{i,\cC}
\times\prod_{\al\in\cC} y_\al
\ee
where, for $\al=(\al_1,\ldots,\al_k)$ in $\cC$, $y_\al$ denotes
$y_{\al_1\ldots\al_k}^{[k]}$.

We now have the following

\begin{lemma}
For any cactus $\cA$, the quantity
$\cY_{i,\cC}$ is independent of the choice of a root $i$ in $[n]$
and of the choice of a refined cactus $\cC$ such that
$\cA=\cA(\cC)$. We will therefore write
$\cY_{\cA}\eqdef\cY_{i,\cC}$ for any such choice
of $i$ and $\cC$.
\end{lemma}

\noindent{\bf Proof :}
First we show the independence with respect to $\cC$.
Let $\cA$ be a cactus, $i$ a fixed root and let
$\cC$ and $\cC'$ be two refined cacti with
$\cA(\cC)=\cA(\cC')=\cA$.
For the given root $i\in [n]$, let $\pi$ be a sequence
constructed as before from the circulations of the
$\al$'s in $\cC$, and let $\pi'$ be an analogous sequence for
$\cC'$.
We need to compare the signatures of $\pi$ and $\pi'$.
For each set $A\in\cA$ of cardinality $k$, there is a unique
$\al\in\cC$ such that ${\til\al}=A$ and a unique $\beta\in \cC'$
such that ${\til\beta}=A$; besides $\beta=(\beta_1,\ldots,\beta_k)$
is a permutation of $\al=(\al_1,\ldots,\al_k)$. The local roots
of $\al$ and $\beta$ coincide and are given by $j=\al_\mu=\beta_\nu$
for $j\in A$ and $1\le \mu,\nu\le k$.
Now note that the signature of the permutation that transforms
the sequence
\[
{\hat\al}=(\al_{\mu+1},\ldots,\al_k,\al_1,\ldots,\al_{\mu-1})
\]
into
\[
{\hat\beta}=(\beta_{\nu+1},\ldots,\beta_k,\beta_1,\ldots,\beta_{\nu-1})
\]
is the same as that which transforms
\[
(j,\al_{\mu+1},\ldots,\al_k,\al_1,\ldots,\al_{\mu-1})
\]
into
\[
(j,\beta_{\nu+1},\ldots,\beta_k,\beta_1,\ldots,\beta_{\nu-1})
\]
Since the latter are respectively circular permutations
of $\al$ and $\beta$, the sign change is the same as the signature of
the permutation that transforms $\al$ into $\beta$.
Indeed a cycle of odd length $k$ has signature $(-1)^{k-1}=1$.
As a result the sign change between the signatures of $\pi$ and
$\pi'$ is exactly compensated by that between
$\prod_{\al\in\cC} y_\al$ and $\prod_{\beta\in\cC'} y_\beta$,
by the antisymmetry of the $y$ tensors.
Therefore $\cY_{i,\cC}=\cY_{i,\cC'}$.

Now we take the same refined cactus $\cC$ with $\cA(\cC)=\cA$ and
compare $\cY_{i,\cC}$ and $\cY_{j,\cC}$ 
for two different choices of {\em global} root :
$i$ and $j$.
Let again $\pi$ be a sequence constructed from the circulations
of the elements in $\cC$ with respect to the root $i$,
and let $\pi'$ be an analogous sequence with respect to the choice
of root $j$.
Note again that there is a unique shortest path in the tree
$G(\cA)$ going from $i$ to $j$. Let $\al_1,\ldots,\al_p$ be
the elements of $\cC$ corresponding to the vertices of $\cA$ that
successively appear along this path.
Let, for each $q$, $1\le q\le p$, ${\hat\al}_q^i$ be the circulation
of the sequence $\al_q$ with respect to the root $i$,
and 
${\hat\al}_q^j$ be the one with respect to the root $j$.
It is easy to see that the signature of the permutation transforming
$\pi$ into $\pi'$ is that of the permutation transforming
the ``reduced'' sequence $\pi_{\rm red}\eqdef
i{\hat\al}_1^i\ldots{\hat\al}_p^i$ into 
$\pi'_{\rm red}\eqdef
j{\hat\al}_p^j\ldots{\hat\al}_1^j$
(we used the obvious notation for the concatenation of words or
sequences). Indeed, choosing $i$ or $j$ as a global root
induces the same local roots for the $\al$'s that are not
on the mentioned path.
By way of example, let us take $i=1$ and $j=18$ for a refined cactus
$\cC$ represented by Figure 2. Then, one would have
\be
\pi_{\rm red} =
(1,2,3,10,11,12,13,14,15,18,19)
\ee
and
\be
\pi'_{\rm red} =
(18,19,15,11,14,12,13,2,10,3,1)
\ee

Now note that, by the tree property of $G(\cA)$, for any $q$,
$1\le q\le p-1$, ${\til\al}_q\cap {\til\al}_{q+1}$
is a singleton whose unique element we denote by $l_q$.
Note also that if one chooses $i$ as a global root,
then the local root of $\al_1$ is $i$, and for any $q$,
$2\le q\le p$, the local root of $\al_q$ is $l_{q-1}$.
On the contrary, if one chooses $j$ as a global root, then
the local root of $\al_p$ is $j$ and for any $q$,
$1\le q\le p-1$, the local root of $\al_q$ is $l_q$.
Remark also that there exist $2p$ (possibly empty)
sequences $u_1,\ldots,u_p$ and $v_1,\ldots,v_p$
such that for any $q$,
$1\le q\le p-1$, ${\hat\al}_q^i$ is equal to the concatenation
$u_q l_q v_q$, while ${\hat\al}_p^i=u_p j v_p$.
One also has ${\hat\al}_q^j=v_q l_{q-1} u_q$, for $2\le q\le p$ and
${\hat\al}_1^j=v_1 i u_1$. 

For the example given by Figure 2, with $i=1$ and $j=18$,
one has
$p=4$, $l_1=2$, $l_2=11$, $l_3=15$,
$u_1=\emptyset$, $v_1=(3)$,
$u_2=(10)$, $v_2=(12,13)$,
$u_3=(14)$, $v_3=\emptyset$,
$u_4=\emptyset$ and $v_4=(19)$.

As a result, we need to evaluate
the signature of the permutation that transforms
\be
\pi_{\rm red}
=
i u_1 l_1 v_1
 u_2 l_2 v_2\ldots
 u_{p-1} l_{p-1} v_{p-1}
 u_p j v_p
\ee
into
\be
\pi'_{\rm red}
=
j v_p l_{p-1} u_p
 v_{p-1} l_{p-2} u_{p-1}\ldots
 v_2 l_1 u_2
 v_1 i u_1
\ee
Notice that one can transform, with a permutation of positive
signature, $\pi_{\rm red}$ into
\be
{\Br\pi}_{\rm red}\eqdef
i u_1 l_1 u_2 l_2\ldots u_{p-1}l_{p-1}
u_p j v_p v_{p-1}\ldots v_1
\ee
This can be done in a succession of steps.
First one changes the segment
$u_1 l_1 v_1 u_2 l_2 v_2$ into $u_1 l_1 u_2 l_2 v_2 v_1$ which
gives a sign $(-1)^{|v_1|(|u_2|+|l_2|+|v_2|)}$
where $|.|$ denotes the length of a sequence.
But $|u_2|+|l_2|+|v_2|=k_2-1$ where $k_2$ is the odd length
of $\al_2$. Then one changes the slightly bigger resulting segment
\[
u_1 l_1 u_2 l_2 v_2 v_1 u_3 l_3 v_3
\]
into
\[
u_1 l_1 u_2 l_2 u_3 l_3 v_3 v_2 v_1
\]
which gives a sign
\be
(-1)^{(|v_2|+|v_1|)(|u_3|+|l_3|+|v_3|)}=1
\ee
since
$|u_3|+|l_3|+|v_3|=k_3-1$ where $k_3$ is the odd length
of $\al_3$ etc.
One can do the same operations with
$\pi'_{\rm red}$ to obtain, without change of sign, the sequence
\be
{\Br\pi}'_{\rm red}\eqdef
j v_p l_{p-1} v_{p-1} l_{p-2}\ldots
v_2 l_1 v_1 i u_1 u_2\ldots u_p
\ee
In the last sequence, one can move $l_1$ in order to lie
between $u_1$ and $u_2$ which gives a factor $(-1)^{|v_1|+|u_1|+1}=1$.
Then we move $l_2$ to make it lie between $u_2$ and $u_3$ which
gives a factor
\be
(-1)^{(|v_1|+|u_1|+1)+(|v_2|+|u_2|+1)}=1
\ee
etc.
Finally, the resulting sequence
\[
j v_p v_{p-1}\ldots v_1
i u_1 l_1 u_2 l_2\ldots u_{p-1} l_{p-1} u_p
\]
can be transformed into ${\Br\pi}_{\rm red}$ by a cycle of length
$n$ and signature $(-1)^{n-1}=1$. This concludes
the proof that $\pi$ is transformed into $\pi'$ by a permutation
of positive signature, and the proof of the lemma.
\endproof

The result of the lemma allows us to state the following
\begin{theor}
The Berezin integral
\[
\int {\rm d}\ch_n\ldots{\rm d}\ch_1\ \ch_i
\exp\lp
\sum_{{3\le k\le n}\atop{k\ {\rm odd}}}
\frac{1}{(k-1)!}
\sum_{(\al_1,\ldots,\al_k)\in[n]^k} y_{\al_1\ldots\al_k}^{[k]}
\ch_{\al_2}\ch_{\al_3}\ldots\ch_{\al_k}
\rp
\]
is independent of $i\in[n]$
and is equal to
\[
\sum_{\cA} \cY_\cA
\]
where the sum is over all odd cacti $\cA$.
\end{theor}

\noindent{\bf Remark 1 :}
In the special case where all the $y$ tensors are zero except for a
specific odd integer $k$, $3\le k\le n$,
one obtains
\be
\int {\rm d}\ch_n\ldots{\rm d}\ch_1\ \ch_i
\exp\lp
\frac{1}{(k-1)!}
\sum_{(\al_1,\ldots,\al_k)\in[n]^k} y_{\al_1\ldots\al_k}
\ch_{\al_2}\ch_{\al_3}\ldots\ch_{\al_k}
\rp
\label{funcint}
\ee
as a sum over all $k$-regular cacti (i.e. cacti $\cA$ whose elements
are subsets of $[n]$ of cardinality $k$). Let the tensor
$A=(A_{\al_2\ldots\al_k})_{(\al_2,\ldots,\al_k)\in[n]^{k-1}}$
be defined by
$A_{\al_2\ldots\al_k}\eqdef \sum_{\al_1=1}^n
y_{\al_1\ldots\al_k}$ and denote by $A^{(i)}$
the tensor obtained from $A$ by forbidding the index $i$.
It is easy to check that (\ref{funcint}) is equal to
\[
(-1)^{i-1} {\rm Pf}^{[k-1]}\lp
A^{(i)}
\rp
\]
where ${\rm Pf}^{[k-1]}(A^{(i)})$ is the order $(k-1)$
{\em Hyperpfaffian} of $A^{(i)}$ as considered, for instance,
in~\cite{Barvinok,Luque}. Note that the result is zero unless $n-1$ is
a multiple of $k-1$.

\noindent{\bf Remark 2 :}
The special case $k=3$ of the previous remark is exactly the
{\em Pfaffian-tree} theorem of Masbaum and Vaintrob~\cite{Masbaum1}.

\noindent{\bf Proof of theorem 2 :}
Our own variation on Brydges' trick is to perform,
for each odd $k$, $3\le k\le n$, and each sequence of indices
$\al_1,\ldots,\al_k$ in $[n]$, the following computation.
Expand
\bea
\lefteqn{(\ch_{\al_1}-\ch_{\al_2})
(\ch_{\al_1}-\ch_{\al_3})\ldots
(\ch_{\al_1}-\ch_{\al_k})=
} & & \nonumber\\
 & & (-1)^{k-1}
\ch_{\al_2}\ch_{\al_3}\ldots\ch_{\al_k}
+(-1)^{k-2}
\sum_{\mu=2}^k
\ch_{\al_2}\ldots\ch_{\al_{\mu-1}}\ch_{\al_1}
\ch_{\al_{\mu+1}}\ldots\ch_{\al_k}
\eea
Notice that for any $\mu$, $2\le \mu\le k$,
\be
(\ch_{\al_2}\ldots\ch_{\al_{\mu-1}})\ch_{\al_1}
(\ch_{\al_{\mu+1}}\ldots\ch_{\al_k})=
\ep_\mu
(\ch_{\al_{\mu+1}}\ldots\ch_{\al_k})
\ch_{\al_1}
(\ch_{\al_2}\ldots\ch_{\al_{\mu-1}})
\ee
with
\be
\ep_\mu=
(-1)^{(\mu-2)(k-\mu)+(\mu-2)+(k-\mu)}=-1
\ee
since $k$ is odd.
As a result, we get
\bea
\lefteqn{
(\ch_{\al_1}-\ch_{\al_2})
(\ch_{\al_1}-\ch_{\al_3})\ldots
(\ch_{\al_1}-\ch_{\al_k})=
} & & \nonumber\\
 & & \sum_{\mu=1}^k
\ch_{\al_{\mu+1}}\ldots\ch_{\al_k}
\ch_{\al_1}\ldots\ch_{\al_{\mu-1}}
\label{idmalek}
\eea
i.e. one obtains, by expanding the product, all the
monomials deduced from $\ch_{\al_2}\ch_{\al_3}\ldots\ch_{\al_k}$
by circular permutation on the {\em full} sequence
$(\al_1,\ldots,\al_k)$.
Since the antisymmetric tensor $y_{\al_1\ldots\al_k}^{[k]}$ is
invariant by circular permutation of its indices ($k$ is odd)
one obtains, writing ${\rm d}\ch$ for
${\rm d}\ch_n\ldots{\rm d}\ch_1$,
\be
\Om_i\eqdef
\int {\rm d}\ch
\ \ch_i
\exp\lp
\sum_{{3\le k\le n}\atop{k\ {\rm odd}}}
\frac{1}{(k-1)!}
\sum_{\al_1,\ldots,\al_k=1}^n y_{\al_1\ldots\al_k}^{[k]}
\ch_{\al_2}\ldots\ch_{\al_k}
\rp
\ee
\bea
\lefteqn{=\int {\rm d}\ch
\ \ch_i
} & & \nonumber\\
 & & \exp\lp
\sum_{{3\le k\le n}\atop{k\ {\rm odd}}}
\frac{1}{k!}
\sum_{\al_1,\ldots,\al_k=1}^n y_{\al_1\ldots\al_k}^{[k]}
(\ch_{\al_1}-\ch_{\al_2})(\ch_{\al_1}-\ch_{\al_3})
\ldots(\ch_{\al_1}-\ch_{\al_k})
\rp \nonumber\\
 & & \,
\eea
\bea
\lefteqn{=\int {\rm d}\ch
\ \ch_i
} & & \nonumber\\
 & & \prod_{{3\le k\le n}\atop{k\ {\rm odd}}}
\lp
\prod_{(\al_1,\ldots,\al_k)\in[n]^k}
\lp
1+\frac{y_{\al_1\ldots\al_k}^{[k]}}{k!}
(\ch_{\al_1}-\ch_{\al_2})(\ch_{\al_1}-\ch_{\al_3})
\ldots(\ch_{\al_1}-\ch_{\al_k})
\rp
\rp \nonumber\\
 & & \,
\eea
Using the antisymmetry of the $y$ tensors,
one can restrict to sequences
$(\al_1,\ldots,\al_k)$ made of distinct elements.
Thus
\be
\Om_i=
\sum_\cC
\lp
\prod_{(\al_1,\ldots,\al_k)\in\cC}
\frac{y_{\al_1\ldots\al_k}^{[k]}}{k!}
\rp
\Om_{i,\cC}
\label{omega}
\ee
where the sum over $\cC$ is over all subsets of
${\cO}_n$ and
\be
\Om_{i,\cC}
\eqdef
\int {\rm d}\ch
\ \ch_i
\lp
\prod_{(\al_1,\ldots,\al_k)\in\cC}
(\ch_{\al_1}-\ch_{\al_2})(\ch_{\al_1}-\ch_{\al_3})
\ldots(\ch_{\al_1}-\ch_{\al_k})
\rp
\label{contrib}
\ee
If two distinct elements $\al$, $\beta$ in $\cC$ are
such that the unordered sets ${\til\al}$ and ${\til\beta}$
coincide then $\Om_{i,\cC}=0$.
Indeed there would then be a permutation $\si$ of $[k]$
for which
$\beta=(\beta_1,\ldots,\beta_k)=(\al_{\si(1)},\ldots,\al_{\si(k)})$ and
the product to be integrated would contain the following
product of $2(k-1)$ factors
\bea
\lefteqn{
(\ch_{\al_1}-\ch_{\al_2})(\ch_{\al_1}-\ch_{\al_3})
\ldots(\ch_{\al_1}-\ch_{\al_k})
} & & \nonumber\\
 & & \times (\ch_{\al_{\si(1)}}-\ch_{\al_{\si(2)}})
(\ch_{\al_{\si(1)}}-\ch_{\al_{\si(3)}})
\ldots
(\ch_{\al_{\si(1)}}-\ch_{\al_{\si(k)}})
\nonumber 
\eea
If $\si(1)=1$ then clearly one can factor, for instance,
$(\ch_{\al_1}-\ch_{\al_2})^2=0$.
If $\si(1)\neq 1$, let $\mu$ be any index, $2\le \mu\le k$, such that
$\si(\mu)\neq 1$ (recall that $k\ge 3$).
One then finds, among the last $k-1$ factors,
$(\ch_{\al_{\si(1)}}-\ch_{\al_{\si(\mu)}})$ which we expand as
\[
(\ch_{\al_1}-\ch_{\al_{\si(\mu)}})-
(\ch_{\al_1}-\ch_{\al_{\si(1)}})
\]
One gets a sum of two terms that vanish since they contain
$(\ch_{\al_1}-\ch_{\al_{\si(\mu)}})^2$ or
$(\ch_{\al_1}-\ch_{\al_{\si(1)}})^2$ which are zero by the
Pauli exclusion principle.
Furthermore, if there is a cycle in $G(\cA(\cC))$
then $\Om_{i,\cC}=0$.
Indeed, there would then be a cycle in the
{\em multigraph} (repeated edges are allowed)
made by putting together all the edges
$\{\al_1,\al_2\}$, $\{\al_1,\al_3\}$,\ldots,
$\{\al_1,\al_k\}$, for each $\al=(\al_1,\ldots,\al_k)$ in $\cC$.
Since, for each edge $\{u,v\}$ in the multigraph, there is
a corresponding factor $(\ch_u-\ch_v)$
in the integrand, the same argument based on telescopic sums
and the Pauli exclusion principle as in the proof of Lemma 1
would show that $\Om_{i,\cC}$ vanishes.
Finally, note that if $G(\cA(\cC))$ does not connect the set $[n]$,
then the integrand of (\ref{contrib}) will not contain
some of the variables $\ch_1,\ldots,\ch_n$, and
$\Om_{i,\cC}$ would be zero.
This shows that the sum in (\ref{omega})
is over refined cacti $\cC$.
For such a $\cC$, one can write, using (\ref{idmalek})
\be
\Om_{i,\cC}=
\int {\rm d}\ch
\ \ch_i
\prod_{(\al_1,\ldots,\al_k)\in\cC}
\lp
\sum_{\mu=1}^k
\ch_{\al_{\mu+1}}\ldots\ch_{\al_k}
\ch_{\al_1}\ldots\ch_{\al_{\mu-1}}
\rp
\ee
One then completely expands the last product, and notices that,
again thanks to the Pauli exclusion principle,
only one term contributes.
Indeed, if the root $i$ belongs to ${\til\al}$, then the only term
in
\[
\sum_{\mu=1}^k
\ch_{\al_{\mu+1}}\ldots\ch_{\al_k}
\ch_{\al_1}\ldots\ch_{\al_{\mu-1}}
\]
which does not contain $\ch_i$ is that for the only index $\mu$
such that $\al_\mu=i$.
We do the same for $\al$'s of ``second generation'' i.e.
which contain an element from a $\beta\in\cC$ such that
$i\in{\til\beta}$, etc.
The end result is that
\be
\Om_{i,\cC}=\int {\rm d}\ch
\ \ch_{\pi_1}\ch_{\pi_2}\ldots\ch_{\pi_n}
\ee
where $\pi=(\pi_1,\ldots,\pi_n)$ is a sequence constructed
using the circulations of elements of the
refined cactus $\cC$, with respect to the choice of root $i$, like
in the considerations preceding the statement of Lemma 2.
Note that $\Om_{i,\cC}=\ep(\pi)$ i.e. the signature of $\pi$
viewed as a permutation of $[n]$. Besides the product
of the $\frac{1}{k!}$ factors in (\ref{omega}) simply accounts
for the number of refined cacti $\cC$ corresponding to the
same cactus $\cA$.
Lemma 2, allows us to write
\be
\Om_i=\sum_{\cA} \cY_\cA
\ee
where the sum
is over all cacti $\cA$, thereby proving the theorem.
\endproof

\end{document}